\documentclass[12pt,a4paper]{article}
\usepackage{amsmath,amssymb,amsfonts}
\def\Bbb{\mathbb}

\title{\bf On a recent reciprocity formula for Dedekind sums}

\author{Kurt Girstmair}
\date{}

\makeatletter
\let\@@maketitle=\maketitle
\def\maketitle{\def\thispagestyle##1{\relax}\@@maketitle}
\makeatother
\newtheorem{theorem}{Theorem}

\newtheorem{corollary}{Corollary}

\def\BE{\begin{equation}}
\def\EE{\end{equation}}
\def\BD{\begin{displaymath}}
\def\ED{\end{displaymath}}
\def\BA{\begin{array}}
\def\EA{\end{array}}
\def\BEA{\begin{eqnarray*}}
\def\EEA{\end{eqnarray*}}
\def\BI{\bibitem}

\def\Z{\Bbb Z}

\def\R{\Bbb R}

\def\phi{\varphi}
\def\EPS{\varepsilon}

\def\MB{\mbox}

\def\MN{\medskip\noindent}

\def\DED{Dedekind }

\begin{document}
\maketitle

\begin{abstract}

\noindent
Let $s(a,b)$ denote the classical Dedekind sum and $S(a,b)=12s(a,b)$. Recently, Du and Zhang proved
the following reciprocity formula. If $a$ and $b$ are odd natural numbers, $(a,b)=1$, then
\BD
  S(2a^*,b)+S(2b^*,a)=\frac{a^2+b^2+4}{2ab}-3,
\ED
where $aa^*\equiv 1\mod b$ and $bb^* \equiv 1 \mod a$. In this paper we show that this formula is a special case
of a series of similar reciprocity formulas. Whereas Du and Zhang worked with the connection of Dedekind sums and values of $L$-series,
our main tool is the three-term relation for Dedekind sums.

\end{abstract}

\section*{1. Introduction and Result}

Let $a$ be an integer, $b$ a natural number, and $(a,b)=1$. The classical \DED sum $s(a,b)$ is defined by
\BD
   s(a,b)=\sum_{k=1}^{b} ((k/b))((ak/b)).
\ED
Here
\BD
  ((x))=\begin{cases}
                 x-\lfloor x\rfloor-1/2 & \MB{ if } x\in\R\smallsetminus \Z; \\
                 0 & \MB{ if } x\in \Z
        \end{cases}
\ED
(see \cite[p. 1]{RaGr}).
It is often more convenient to work with
\BD
 S(a,b)=12s(a,b)
\ED instead. We call $S(a,b)$ a {\em normalized} \DED sum.

Probably the most important elementary result concerning \DED sums is reciprocity law. If $a$ and $b$ are coprime natural numbers, then
\BE
\label{1.2}
   S(a,b)+S(b,a)=\frac{a^2+b^2+1}{ab}-3.
\EE
Recently, Du and Zhang have found the following hitherto unknown reciprocity law (see \cite{DuZh}).
If $a$ and $b$ are coprime odd natural numbers, then
\BE
\label{1.4}
  S(2a^*,b)+S(2b^*,a)=\frac{a^2+b^2+4}{2ab}-3,
\EE
where $aa^*\equiv 1\mod b$ and $bb^* \equiv 1 \mod a$.

The proof given in \cite{DuZh} is based on the connection of \DED sums and values of $L$-series. The authors of the said paper ask for an elementary proof of their
result. Here we give such an elementary proof based on the tree-term-relation of \DED sums. Moreover, we show that (\ref{1.4}) is a special case of a series of similar
reciprocity formulas. Indeed, we have the following.

\begin{theorem} 
\label{t1}

Let $a$ and $b$ be coprime natural numbers and $t$ a natural number such that $a^2+1\equiv 0\mod t$. Further, let $(b,t)=1$.
Then
\BE
\label{1.6}
 S(ta^*,b)+S(tb^*,a)=\frac{a^2+b^2+t^2}{tab}-3+S(ab,t).
\EE
\end{theorem} 

\MN
As to the case $t=1$, we note
\BE
\label{1.7}
S(a^*,b)=S(a,b)
\EE
(see \cite[p. 26]{RaGr}) and $S(ab,1)=0$. In the case $t=2$, $a$ and $b$ are odd and $S(ab,2)=0$. Hence we obtain the
following.

\begin{corollary} 
\label{c1}
The formulas {\rm(\ref{1.2})} and {\rm(\ref{1.4})} are immediate consequences of Theorem \ref{t1} in the cases $t=1$ and $t=2$.
\end{corollary} 

\begin{corollary} 
\label{c2}
Suppose, in the setting of Theorem \ref{t1}, that $b\equiv\pm 1\mod t$. Then
\BE
\label{1.8}
 S(ta^*,b)+S(tb^*,a)=\frac{a^2+b^2+t^2}{tab}-3.
\EE
Suppose, on the other hand, that $b\equiv\pm a\mod t$. Then
\BE
\label{1.10}
 S(ta^*,b)+S(tb^*,a)=\frac{a^2+b^2+t^2}{tab}+\begin{cases} -t-2/t ,  & \MB{ if } b\equiv a \mod t;\\
                                                            t+2/t-6, & \MB{ if } b\equiv -a\mod t.
                                                     \end{cases}
\EE
\end{corollary} 

\MN
As to (\ref{1.8}), note that $(ab)^2\equiv-1\mod t$, which shows that $S(ab,t)=0$ (see \cite[p. 28]{RaGr}).
In the case of (\ref{1.10}), we use $S(1,t)=t+2/t-3$ (which is an immediate consequence of (\ref{1.2})) and $S(-a,b)=-S(a,b)$ (see \cite[p. 26]{RaGr}).

\MN
{\em Remark.} The natural numbers $t$ such that there is a natural number $a$ with $a^2+1\equiv 0 \mod t$ can be characterized as follows: $t=m$ or $t=2m$, where $m$ is a natural number
whose prime divisors are all $\equiv 1 \mod 4$ (this includes $m=1$).

\MN
{\em Example}. Let $t=5$, $a$ such that $a^2+1\equiv 0\mod 5$, and $(a,b)=(b,5)=1$. This implies $a\equiv \pm 2\mod 5$. If $b\equiv \pm 1\mod 5$, then
\BD
   S(5a^*,b)+S(5b^*,a)=\frac{a^2+b^2+25}{5ab}-3.
\ED
In the remaining case, we have $b\equiv \pm a\mod 5$. If $b\equiv a\mod 5$, then (\ref{1.10}) reads
\BD
   S(5a^*,b)+S(5b^*,a)=\frac{a^2+b^2+25}{5ab}-27/5.
\ED
If $b\equiv-a\mod 5$, we have
\BD
   S(5a^*,b)+S(5b^*,a)=\frac{a^2+b^2+25}{5ab}-3/5.
\ED

\section*{Proof of Theorem \ref{t1}}

Let $a,b,t$ be natural numbers, $(a,b)=(b,t)=1$, such that $a^2+1\equiv 0\mod t$. Put $c=b(a^2+1)/t$. Obviously, $(a,c)=1$. Then \cite[Th.4 ]{Gi1} says
\BE
\label{2.2}
    S(a,c)=\frac{(b^2-1)a}{tb}-S(ab,t)+S(at^*,b).
\EE
where $tt^*\equiv 1\mod b$.
By the reciprocity law (\ref{1.2}),
\BD
   S(a,c)=-S(c,a)+\frac{a^2+c^2+1}{ac}-3.
\ED
However, $c\equiv bt^*\mod a$, with $tt^*\equiv 1 \mod a$. Hence $S(c,a)=S(bt^*,a)$. We replace $at^*$ by $ta^*$ and $bt^*$ by $tb^*$
in the respective normalized \DED sums (see (\ref{1.7})).
Then a short calculation proves Theorem \ref{t1}.

We still have to make clear that this remarkably simple proof is based on elementary results. Indeed, (\ref{2.2}) follows from the three-term relation
\BD
  S(a,b)=S(c,d)+\EPS S(r,|q|) +\frac{b^2+d^2+q^2}{bdq}-3\EPS
\ED
(see \cite{Gi1}).
Here $b,d,$ are natural numbers, $a,c$ integers, $(a,b)=(c,d)=1$, $a/b\ne c/d$. Further, $q=ad-bc$ and $\EPS$ is the sign of $q$. Finally,
$r=aj-bk$, where $j,k$ are integers such that $-cj+dk=1$. The three-term relation, in turn, can be deduced from the composition rule of
the logarithm of Dedekind's $\eta$-function (see \cite{Di, Gi2}). An elementary proof of this composition rule is given in \cite[\S4]{Ra}.


\vspace{0.5cm}
\noindent
Kurt Girstmair            \\
Institut f\"ur Mathematik \\
Universit\"at Innsbruck   \\
Technikerstr. 13/7        \\
A-6020 Innsbruck, Austria \\
Kurt.Girstmair@uibk.ac.at

\end{document}